\documentclass[11pt]{amsart}
\usepackage{amssymb}
\begin{document}
\title{Remarks about Schlumprecht space}
\author{Denka Kutzarova}
\address {Institute of Mathematics\\ Bulgarian Academy of Sciences\\ 1113
Sofia, Bulgaria}
\curraddr{Miami University, Oxford, Ohio, USA\\e-mail address:
kutzard@muohio.edu}

\begin{abstract}
We prove that Schlumprecht space $S$ is isomorphic to
$\left(\sum^\infty_{k=1} \oplus\ \ell_\infty^{n_k}\right)_S$ for any
sequence of integers $\left(n_k\right)$.
 \end{abstract}
\maketitle
Th. Schlumprecht constructed the first arbitrarily distortable space
$S$ \cite{S1}.  He also proved that $S$ is complementably minimal \cite{S2}.  A
Banach space $X$ is called minimal (a notion defined by H.P. Rosenthal)
if every infinite dimensional subspace of $X$ contains a further subspace
isomorphic to $X$ and $X$ is called complementably minimal (defined by A.
Pe\l czynski) if every infinite dimensional subspace of $X$ contains a
subspace which is isomorphic to $X$ and complemented in $X$.

As it was remarked in \cite{S2}, the space $S$ is either prime or fails the
Schroeder-Bernstein property (see \cite{C}, \cite{G2}).  A Banach space $X$ is called
prime \cite{LT} if every complemented infinite dimensional subspace of $X$ is
isomorphic to $X$.  A. Pe\l czynski showed that the spaces $\ell_p, 1 \leq
p < \infty$, and $c_0$ are prime, while J. Lindenstrauss proved that
$\ell_\infty$ is also prime (cf. e.g. \cite{LT}).  Recently W.T. Gowers and B.
Maurey \cite{GM} constructed more prime spaces for which all complemented
infinite dimensional subspaces are of finite codimension.  It is not known
yet if $S$ is prime.  It seems hard to characterize its complemented
subspaces since $S$ is rich with them.  We shall see that for every
sequence $\left(n_k\right)^{\ \infty}_{k=1}$ of integers, the Schlumprecht
sum $Y = \left(\sum\limits^\infty_{k=1} \oplus\ \ell_\infty^{n_k}\right)_S$
is isomorphic to a complemented subspace of $S$.  For some sequences, one
can immediately observe that the corresponding sum is isomorphic to $S$.
For example, if $n_k = 2^k$, we use that the standard basis of $S$ is
subsymmetric to show that $Y$ is isomorphic to its square and then it
follows by Pe\l czynski decomposition technique that $Y \approx S$ (see
\cite{C} for more details).  We hoped in the beginning that there exists a sequence
$\{n_k\}$ for which $Y$ is not isomorphic to $S$ and thus to show that $S$
fails to be prime.  Unfortunately, this turned out not to be true.

Recall the definition of Schlumprecht space $S$. Let
$\left(e_i\right)^{\ \infty}_{i=1}$ be the standard basis of the linear
space $c_{00}$ of finitely supported sequences.  For $x =
\sum\limits^\infty_{i=1} a_i e_i \in c_{00}$, set supp $x = \{i \in
\mathbb{N} : a_i \not= 0\}$.  For $
E,F$ finite subsets of $\mathbb{N},\ E < F$
means max $E < \min F$ or either $E$ or $F$ is empty.  For
$x = \sum\limits^\infty_{i=1} a_ie_i$ and $E$ a subset of $\mathbb{N}$, we
denote by $Ex$ the vector $Ex = \sum\limits_{i \in E} a_ie_i$.  For $t \geq
1$ let $f(t) = \log_2(t+1)$.  Then $S = (S, \Vert\cdot\Vert)$ is the
completion of $c_{00}$ with respect to the norm $\Vert \cdot\Vert$ defined
by the equation:

\begin{displaymath}
\Vert x \Vert = \max\left\{\Vert x \Vert_\infty,
\sup_{\substack
{r \geq 2\\
E_1<E_2<\cdots< E_r}}
\ \ \  \frac 1{f(r)} \sum^r_{i=1} \Vert
E_ix\Vert \right\}\ .
\end{displaymath}

Clearly, $(e_i)$ is a 1-unconditional and 1-subsymmetric basis of $S$.

We shall use the following results from \cite{S2}.

{\bf{Statement 1.}}  Every normalized block basis $(y_i)$ of $(e_i)$
dominates $(e_i)$, i.e.
\begin{displaymath}
\Vert \sum^\infty_{i=1}a_iy_i \Vert \geq \Vert \sum^\infty_{i=1}
a_i e_i\Vert
\end{displaymath}
for all $(a_i) \in c_{00}$.

\smallskip
{\bf{Statement 2.}} (Lemma 2 \cite{S2})  Let $\varepsilon > 0$ and $r \in
\mathbb{N}$.  Then there is an $n = n(\varepsilon,r) \in \mathbb {N}$ with
following property: if $m \geq n$ and if $y = \frac 1m \sum\limits^m_{i=1}
x_i$, where $(x_i)^{\ m}_{i=1}$ is a normalized block basis of $(e_i)$
which is
$(1 + \varepsilon/2)$ - equivalent to the unit basis of $\ell^m_1$, then
\begin{displaymath}
\sup_{E_1<E_2<\dots<E_r} \ \ \sum^r_{i=1} \Vert E_iy\Vert  \leq \Vert y \Vert +
\varepsilon .
\end{displaymath}

\smallskip
{\bf{Statement 3.}}  (Theorem 3 \cite{S2})  Let $\varepsilon_k > 0, k=1,2,\dots$
with $\Sigma\ \varepsilon_k \leq 1$.  There exists a constant $C > 1$ such
that for any normalized block basis $(y_k)$ of $(e_k)$ with the
following properties:  there is a sequence $r_k \uparrow \infty$ in $\mathbb N$ so
that for all $k \in \mathbb N$,

\smallskip
(a)\qquad $\displaystyle\sup_{\substack{r \leq r_{k-1}\\ E_1 < E_2 <\dots<E_r}}
\sum^r_{i=1} \Vert E_i(y_k) \Vert \leq 1 + \varepsilon_k$\ ,

\smallskip
(b)\qquad $\displaystyle \text{card} (\text{supp } y_k) \leq \varepsilon_k \cdot
f\left(\frac {r_k}3\right)$,

\smallskip
\noindent
we have that $(y_k)$ is $C$-equivalent to $(e_k)$.

The above statements easily imply that $S$ is complementably minimal.  It
is known that $S$ contains uniformly $\ell^n_\infty, n = 1,2,\dots$  The
construction of a $(1 + \delta)$-copy of $\ell^n_\infty$ looks like a
``yardstick''.  For example, let
$F_1, F_2,\dots,F_n$ be disjoint subsets of $\mathbb N$ with card$(F_1) =
p_1 = q_1$, card $(F_2) = p_2 = q_1q_2,\dots,$ card $(F_n) = p_n =
q_1 q_2 \dots q_n$.  We also construct $F_1,\dots, F_n$
to be ``nested'', i.e. for $\chi_{F_i} = \{e^i_{k_j}\}^{p_i}_{j=1}\ , i =
1,2,\dots,n,$ we require supp$\{e_{k_1}^i\} < \text{supp}\{e_{k_j}^{i+1} :
j = 1,2,\dots,q_{i+1}\} < \text{supp} \{e_{k_2}^i\} < \text{supp}
\{e_{k_j}^{i+1} : j = q_{i+1} + 1,\dots, 2 q_{i+1}\} < \dots <
\text{supp}\{e^i_{k_{p_i}}\} < \text{supp}\{e_{k_j}^{i+1} : j = (p_i - 1)q_{i+1} ,\dots, p_{i+1}\}\ .$

If $q_i$ are chosen big enough (depending on $n$ and $\delta$), $q_1 \ll
\dots \ll q_n$, then the finite sequence of normalized vectors
$y_i = \frac 1{f(p_i)}\chi_{F_i},\ i = 1,2,\dots,n$, is
$(1+\delta)$-equivalent to the unit basis of $\ell_\infty^n$.  A proof of
this fact was never written.  One can find an argument in the same spirit
in \cite{ADKM}; it is more complicated, however, because mixed Tsirelson
spaces are asymptotic $\ell_1$ spaces.

Let now $(n_k)_{k=1}^{\ \infty}$ be an arbitrary sequence of integers.  We
shall construct a subspace $Y$ of $S$ which is isomorphic to $(\Sigma\
\oplus \ell_\infty^{n_k})_S$.  Recall that given $(X_k, \Vert \cdot
\Vert_k)$, the norm of $y = \Sigma\ y_k$ in the sum is defined as $\Vert
\Sigma \Vert y_k \Vert_k e_k\Vert$, where $\Vert \cdot \Vert$ is the
standard norm of $S$.

In view of Statement 3, choose a sequence of positive reals $\varepsilon_k$
with $\Sigma\ \varepsilon_k \leq 1$.  We shall construct successively
for $k = 1,2,\dots$ collections $(z_i^k)_{i=1}^{n_k}$ of normalized
elements in $S$ with $[(z^k_i)_{i=1}^{n_k}]$ being $(1 + \delta_k)$-isomorphic to
$\ell_\infty^{n_k}$ for sufficiently small $\delta_k$ such that $y_k =
z_k/\Vert z_k\Vert$, where $z_k = \sum\limits^{n_k}_{i=1} z^k_i$, satisfy
the conditions of Statement 3.  Each $z_k$ will be an arithmetic mean of
$m_k$ ``yardsticks'', $m_k$ big enough and $m_k \ll q^k_1 \ll q_2^k \ll\dots\ll q^k_{n_k}$.  That
is, we pick a yardstick, described above, then we take $m_k$ successive
copies of it, i.e. $m_k$ successive vectors with the same distribution.
Thus, each
$z_i^k, i = 1,\dots, n_k$, is of the form from Statement 2 with $m = m_k$.
Note that if $n_k = n_i$ for $k \not= i$, in order to satisfy the
conditions of Statement 3, the concrete representations of
$\ell_\infty^{n_k}$ and $\ell_\infty^{n_i}$ in $S$ have different
distributions.

By Statement 3, $(y_k)$ is $C$-isomorphic to $(e_k)$.  The form of $y_k$
also assures that if for each $k$ we take an arbitrary combination $x_k =
\sum\limits^{n_k}_{i=1} a_i^k y_i^k$, then the sequence
$(x_k/\Vert x_k\Vert)$ also satisfies the hypothesis of Statement 3.
Therefore,
$(x_k/\Vert x_k\Vert)$ is $C$-isomorphic to $(e_k)$.  Set $Y =
[(y_k)_{k=1}^{\ \infty}]$.  Thus, we obtained

\medskip
{\bf{Proposition 1.}}
The subspace $Y$ is $2C$-isomorphic to $(\sum\limits^\infty_{k=1} \oplus\
\ell_\infty^{n_k})_S$.

We also have

\medskip
{\bf{Proposition 2.}}
The subspace $Y$ is complemented in $S$.

\begin{proof}
Let $E_k$ be intervals of integers, containing supp $y_k$ and $E_1 < E_2
<\dots$.  Since the subspaces $\ell_\infty^n$ are uniformly complemented,
see \cite{LT}, let $Q_k$ be projections from $E_kS$ onto $[(y_i)_{i=1}^{n_k}]$
with $\Vert Q_k\Vert \leq 3$.  It follows from Statement 1 that for every
$x \in S$,
\begin{displaymath}
\Vert \sum^\infty_{k=1} \Vert E_k x \Vert e_k \Vert \leq \Vert
\sum^\infty_{k=1} E_kx \Vert \leq \Vert x \Vert\ .
\end{displaymath}
On the other hand, for the map $P = \sum\limits_k Q_kE_k$ we have by
Proposition 1 that
\begin{displaymath}
\Vert Px\Vert \leq 2C \Vert \sum^\infty_{k=1} \Vert Q_k E_k x \Vert e_k
\Vert \leq 6C \Vert \sum^\infty_{k=1} \Vert E_k x \Vert e_k\Vert\ .
\end{displaymath}
Therefore, $\Vert Px\Vert \leq 6C \Vert x\Vert$.  Clearly, $P^2 = P$ and
$P(S) = Y$, which ends the proof.
\end{proof}

In order to show that $Y$ is isomorphic to $S$, we shall use a variant of
the refinement of Johnson's argument, presented in \cite{BCLT}, Proposition 6.3,
for a space with a subsymmetric basis.  We shall work with finite dimensional
decompositions instead of block-vectors.

\medskip
{\bf{Proposition 3.}}
The subspace $Y$ is isomorphic to $S$.

\begin{proof}
Following \cite{BCLT}, represent every integer in the form $m = 2^i(2j-1),\ i =
0,1,2,\dots, j = 1,2,\dots$\ .  As above, we construct inductively a
sequence of consecutive finite dimensional subspaces $(V_m)$,
i.e. supp $V_1 < \text{supp } V_2 <\dots$, where we define the support of a subspace
to be the union of the supports of its elements.  We construct $(V_m)$
to be generated by an appropriate yardstick which represents in $S$ the
space $\ell^{n_j}_\infty$ and we also assure that the conditions of
Statement 3 are satisfied.  Note that the representation of
$\ell_\infty^{n_j}$ have different distribution for different $m$.  As in
Proposition 1, the subspace $V = [(V_m)_{m=1}^{\ \infty}]$ is isomorphic to
the sum $Z = (\ell_\infty^{n_1} \oplus \ell_\infty^{n_1} \oplus
\ell_\infty^{n_2} \oplus \ell_\infty^{n_1} \oplus \ell_\infty^{n_3} \oplus
\dots)_S$.  By Proposition 2, $V$ is complemented in $S$.  We have that
$[(V_{2m})_{m=1}^{\ \infty}]$ is also isomorphic to $Z$ and therefore, it
is also isomorphic to $V$.  The subspace $[(V_{2k-1})_{k=1}^{\ \infty}]$ is
isomorphic to $(\sum\limits^\infty_{k=1} \oplus\ \ell_\infty^{n_k})_S$ and
hence it is isomorphic to $Y$.  Since $V$ is complemented, $S \approx V
\oplus W$, whence
\begin{displaymath}
Y \oplus S \approx Y \oplus V \oplus W \approx V \oplus W \approx S.
\end{displaymath}
Since $S$ is complementably minimal, $Y \approx S \oplus X$.  Evidently $S
\oplus S \approx S$ and therefore,
\begin{displaymath}
Y \approx S \oplus X \approx S \oplus S \oplus X \approx S \oplus Y \approx
S.
\end{displaymath}
Thus, $Y \approx S$, which implies $S \approx (\Sigma \oplus
\ell_\infty^{n_k})_S$.
\end{proof}

A Banach space with an unconditional basis is said to have a unique
unconditional basis if any two normalized unconditional bases are
equivalent after a permutation, see \cite{BCLT} and for recent results see
\cite{CK}.  It was remarked in \cite{CK} that the example of Gowers in
\cite{G1} has a unique
unconditional basis.  Proposition 3 implies that $S$ has no unique
unconditional basis.  This also immediately follows if one uses that $S$ is
complementably minimal and combines Proposition 6.3 and Proposition 6.4
\cite{BCLT}.  They provide the existence of a normalized block basis with
constant coefficients $(u_i)$ which spans a subspace isomorphic to $S$ and
no permutation of $(u_i)$ is equivalent to $(e_i)$.  In fact, one can
construct such a block-basis $(u_i)$ directly without Proposition 6.4 which
uses an adaptation of Zippin's argument.  Indeed, if we consider consecutive
blocks with constant
coefficients whose length increases rapidly enough, it is easy to show that
Statement 3 and Proposition 6.3 \cite{BCLT} imply the following.

\medskip
{\bf{Proposition 4.}}
Schlumprecht space $S$ is isomorphic to
$(\sum\limits_k \oplus\ \ell_1^{n_k})_S$
for any sequence of integers $(n_k)$.

\medskip
P. Casazza pointed out the following immediate consequence of our observations
(in fact, not using their full generality). Recall that the same property
for the classical spaces $\ell_p$ for $p>2$, $1<p<2$ and $p=1$ was shown in
\cite{R}, \cite{BDGJN} and \cite{B} respectively.

\medskip
{\bf{Corollary 5.}}
Schlumprecht space $S$ (resp. its dual $S^*$) has a subspace isomorphic to
the whole space and not complemented in $S$ (resp. $S^*$).

\medskip
One possible way to see it is to use the fact that one can embed $\ell_1^n$
in $\ell_\infty^m$ with a bad constant of complementation.

\medskip
{\bf{Remark.}}
Let $Z$ be a complemented subspace of $S$ which has a subsymmetric basis.
Since $Z$ is isomorhic to its square, then by the decomposition technique
we have that Z is isomorphic to $S$.

\medskip
{\bf{Question.}}  Does $S$ have a unique subsymmetric basis?
\medskip

\bibliographystyle{plain}

\end{document}